\documentclass{amsart}
\usepackage{amssymb}
\usepackage{amsfonts}
\usepackage{amsmath}
\usepackage{lastpage}
\usepackage[legalpaper,bookmarks=true,colorlinks=true,linkcolor=blue,citecolor=blue]{hyperref}
\usepackage{graphicx}%
\usepackage{fancyhdr}
\usepackage{color}
\usepackage[mathlines]{lineno}
\usepackage{lscape}
\usepackage{epsfig}
\usepackage{geometry}
\usepackage{tgbonum}
\fontfamily{qcr}\selectfont

\newtheorem{theorem}{Theorem}
\theoremstyle{plain}

\newtheorem{corollary}{Corollary}

\newtheorem{proposition}{Proposition}

\numberwithin{equation}{section}

\begin{document}
\Large
\title[Alternative Asymptotic Theory for atom endpoints]{Finiteness of Record values and Alternative Asymptotic Theory of Records with Atom Endpoints}

\begin{abstract} Asymptotic theories on record values and times, including central limit theorems, make sense only if the sequence of records values (and of record times) is infinite. If not, such theories could not even be an option. In this paper, we give necessary and/or sufficient conditions for the finiteness of the number of records. We prove,  for example for \textsl{iid} real valued random variable, that strong upper record values are finite if and only if the upper endpoint is finite and is an atom of the common cumulative distribution function. The only asymptotic study left to us concerns the infinite sequence of hitting times of that upper endpoints, which by the way, is the sequence of weak record times. The asymptotic characterizations are made using negative binomial random variables and the dimensional multinomial random variables. Asymptotic comparison in terms of consistency bounds and confidence intervals on the different sequences of hitting times are provide. The example of a binomial random variable is given.\\

\bigskip
\noindent $^{\dag}$ Gane Samb Lo.\\
LERSTAD, Gaston Berger University, Saint-Louis, S\'en\'egal (main affiliation).\newline
LSTA, Pierre and Marie Curie University, Paris VI, France.\newline
AUST - African University of Sciences and Technology, Abuja, Nigeria\\
gane-samb.lo@edu.ugb.sn, gslo@aust.edu.ng, ganesamblo@ganesamblo.net\\
Permanent address : 1178 Evanston Dr NW T3P 0J9,Calgary, Alberta, Canada.\\

\noindent \noindent $^{\dag\dag}$ Harouna Sangar\'e.\\
DER MI, Facult\'e des Sciences et Techniques, USTTB.\newline
harounasangare@fst-usttb-edu.ml, harouna.sangare@mesrs.ml\\

\noindent \noindent $^{\dag\dag\dag}$ Mamadou Cherif Moctar Traor\'e.\\
LERSTAD, Gaston Berger University, Saint-Louis, S\'en\'egal.\newline
mamdou traore.cherif-mamadou-moctar@ugb.edu.sn.\\

\bigskip
\noindent $^{\dag \dag \dag\dag}$ Mohammad Ahsanullah\\
Department of Management Sciences, Rider University, Lawrenceville, New Jersey, USA\\
ahsan@rider.edu\\

\end{abstract}
\maketitle

\section{Introduction}

\noindent A considerable number of asymptotic results are available in the literature concerning infinite sequences of record values of sequence of real-valued random variables defined on the same probability space. However if the upper-endpoint $x_0$ of the common cumulative distribution function (\textit{cdf}) of a sequence of independent of identically distributed (\textit{iid}) copies of $X$  is an atom of the common probability law $\mathbb{P}_X$, that is $\mathbb{P}_X(\{x_0\})$ is not zero, there will not be any further record value once $x_0$ is hit for the first time. In such situations, all the available asymptotic theories become irrelevant. The paper shows the facts already described above and proposes asymptotic results on the sequences of hitting times of the upper end-point or the lower endpoints.\\

\noindent Let $X_{1}, X_{2}, \cdots$ be a sequence of real-valued random variables defined on the same probability space $\left( \Omega, \mathcal{A}, \mathbb{P}\right)$. Let us defined the sequence of record values:\\
 
$X^{(1)}=X_{1}$, $U(1)=1$, and upon the existence of $Y^{(n)}=X_{U(n)}$
 
$$
X^{(n+1)}= \left\{
\begin{array}{c}
\inf\{j > U(n), X_{j} > Y^{(n)}  \},  \\
\text{Not defined }  \text{ if } \inf \{ j > U(n), X_{j} > Y^{(n)}  \}=+\infty.  
\end{array}
\right.
$$
 
\bigskip \noindent The sequence $\{ U(n), 1 \leq n \leq M  \}$ is that of the occurence times of the record values, where $M$ can be infinite or finite, constant or random.\\
 
\noindent  If there is no further record value after $X^{(n)}$, we denote $U(n+1)=+\infty$.\\
 
\noindent In this note, we begin by a general result on the finiteness or infiniteness of the total number of records regardless the dependence structure of $\left( X_{n} \right)_{n \geq 1}$ in Section \ref{finitenessTNR}. After the justification of the finiteness of upper record if $uep(F)$ is an atom and the finiteness of lower records if $lep(F)$ is an atom, we study the infinite sequences of hitting times of one of the endpoints in Section \ref{sec2}. We  provide central limit theorems, laws of the iterated logarithm and Berry-Essence bounds for the sequence of hitting times $(N_{i,k})_{k\geq 1}$ of $lep(F)$ ($i=1$) and of $uep(F)$ ($i=2$) (Theorem \ref{theo02}). These results are based on probability laws of negative binomial random variable. \\

\noindent Next, in the same section, we take a multinomial approach with three outcomes $E_1$ (lep(F) is hit) of probability $p_1$, $E_2$ (uep(F) is hit) of probability $p_2$ and $E_3$ (neither lep(F) nor uep(F) is hit) of probability $p_3=1-p_1-p_2$. We study the asymptotic law of the mutinomial random vector $(M_{1,n}, \ M_{2,n}, \ M_{3,n})^{t}$ (Proposition \ref{rec-02}), and draw asymptotic sub-results on the difference $M_{1,n}-M_{2,n}$ and on the ratio $M_{1,n}/M_{2,n}$ in Corollary \ref{rec-03}.\\

\noindent In summary, for $p_1-p_2>0$ for example, the difference $M_{1,n}-M_{2,n}$ goes to infinity with the rate $\delta n$, where $0<\delta<p_1-p_2$. In section \ref{sec3}, we give an illustration in the case where $X$ is a binomial random variable $\beta(r,\alpha)$, $r\geq 1$, $\alpha \in ]0,1[$ and provide confidence interval for
$M_{1,n}/M_{2,n}$ in the case where $\alpha=1/2$. We finish the paper by a conclusive section.

\section{Finiteness or Infiniteness of the total number of records} \label{finitenessTNR}

\noindent Let us begin by a general law.

\begin{proposition} \label{nrec01}
For each $k\geq 1$, set 
$$
X^{\star}_{k}= \sup_{h>k} X_h.
$$

\bigskip
\noindent and denote

$$
D_{-}=\{(x,y) \in \mathbb{R}^2, \  x\leq y \}.
$$

\bigskip
\noindent We have for any $n\geq 2$

$$
\mathbb{P}(U(n+1)=+\infty)=\sum_{k\geq n} \mathbb{P}_{(X^{\star}_{k},X_k)}(D_{-}) \mathbb{P}(U(n)=k).
$$
\end{proposition}

\bigskip
\noindent \textbf{Proof}. Conditioning on $(U(n)=k)$, $(U(n+1)=+\infty)$ means that all the $X_h$, $h>k$, are less than $X_k$. The proof is ended by the remark

$$
\mathbb{P}(\max_{h>k} X_h\leq X_k)=\mathbb{P}(X^{\star}_{k}\leq X_k)=\mathbb{P}((X^{\star}_{k},X_k) \in D_{-})=\mathbb{P}_{(X^{\star}_{k},X_k)}(D_{-}). \blacksquare
$$

\noindent Let us give an application of Proposition \ref{nrec01} in the independent case.

\begin{proposition} \label{nrec02}
Suppose that $X_1$, $X_2$, \ldots are independent random variables with respective cumulative distribution functions (\textit{cdf}) $F_{j}$, $j\geq 1$. Then, whenever $U(n)$ is finite, we have

$$
\mathbb{P}(U(n+1)=+\infty)=\sum_{k\geq n} \left( \int_{\mathbb{R}} \biggr( \prod_{j>k} F_{j}(x) \biggr) d\mathbb{P}_{X_k}(x) \right) \mathbb{P}(U(n)=k).
$$
\end{proposition}

\noindent \textbf{Proof}. Here $X^{\star}_{k}$ and $X_k$ are independent and we have

$$
\mathbb{P}(U(n+1)=+\infty)=\sum_{k\geq n} \mathbb{P}_{X^{\star}_{k}} \otimes \mathbb{P}_{X_k}(D_{-}) \mathbb{P}(U(n)=k).
$$

\noindent Let us use Fubini's theorem to have

\begin{eqnarray*}
&&\mathbb{P}_{X^{\star}_{k}} \otimes \mathbb{P}_{X_k}(D_{-})\\
&&=\int_{\mathbb{R}} d\mathbb{P}_{X_k}(x) \int_{\mathbb{R}} 1_{D_{-}}(x,y) d\mathbb{P}_{X^{\star}_{k}}\\
&&=\int_{\mathbb{R}} \mathbb{P}(X^{\star}_{k} \leq x) d\mathbb{P}_{X_k}(x)\\
&&=\int_{\mathbb{R}} \biggr( \prod_{j>k} F_{j}(x) \biggr) d\mathbb{P}_{X_k}(x).
\end{eqnarray*}

\bigskip
\noindent We get the announced result by combining the above lines.\\

\bigskip \noindent Now let us see what happens if the sequence is stationary, that is $F_j=F$ for all $j\geq 1$. Define the lower and the upper endpoints (lep and uep) of $F$ by

$$
lep(F)=\inf \{x \in \mathbb{R}, \ F(x)>1 \} \  and \ uep(F)=\sup \{x \in \mathbb{R}, \ F(x)<1 \}
$$

\bigskip \noindent  We have

$$
\int_{\mathbb{R}} \biggr( \prod_{j>k} F_{j}(x) \biggr) d\mathbb{P}_{X_k}(x)=\int_{-\infty}^{uep(F)} F(x)^{+\infty} dF(x).
$$

\bigskip \noindent \noindent But $F(x)^{+\infty}=0$ unless $x=uep(F)$. This gives

$$
\int_{\mathbb{R}} \biggr( \prod_{j>k} F_{j}(x) \biggr) d\mathbb{P}_{X_k}(x)=\int_{-\infty}^{uep(F)} 1_{\{uep(F)\}} dF(x)=\mathbb{P}(X=uep(F)).
$$

\bigskip \noindent  We conclude that

$$
\mathbb{P}(U(n+1)=+\infty)=\sum_{k\geq n} \mathbb{P}(X=uep(F)) \mathbb{P}(U(n)=k)=\mathbb{P}(X=uep(F))
$$

\noindent which leads to the simple result:

\begin{proposition} \label{nrec03} 

Suppose that $X_1$, $X_2$, \ldots are independent and identically distributed random variables with common \textit{cdf} $F$ and let $uep(F)$ denote the upper endpoint of $F$. Then

$$
\mathbb{P}(U(n+1)=+\infty)=\mathbb{P}(X=uep(F).
$$
\noindent As a consequence, the sequence of record values (and of record times) is finite if and only if $uep(F)$ is finite and is an atom of $F$, 
that is $\mathbb{P}_X(uep(F))>0$.
\end{proposition}

\bigskip \noindent \textbf{Consequences}. The number of time records \textit{a.s.} is infinite in the following cases.\\

\noindent (1) $uep(F)=+\infty$.\\

\noindent (2) $uep(F)<+\infty$ but $\mathbb{P}(X=uep(F))=0$. Example : $X \sim \mathcal{U}(0,1)$.\\

\bigskip \noindent The number of time records may be finite in the following cases.\\

\noindent (1) $X$ is discrete and takes a finite number of points.\\

\noindent (2) $X$ is discrete, takes an infinite number of values such that the strict values set $\mathcal{V}_X$ of $X$ has a maximum value. We mean by strict values set, the set of points taken by $X$ with a non-zero probability.\\

\section{Infinite number of hitting times for the extreme endpoints} \label{sec2}

\noindent  In this note, we focus on $X$, $X_{1}$, $X_{2}$, \ldots are \textit{iid} random variables with \textit{cdf} $F$ such that $uep\left(F\right)$ is finite and is an atom of $F$, that is

\begin{equation}\label{H1a}
\mathbb{P}\left(X=uep\left(F\right)\right)=p_2 \in \left] 0, 1\right[.
\end{equation}

\bigskip \noindent In that context, let us see right now that $uep(F)$ will be hit infinitely many times, that is

$$
\{j\geq 1, \ X_j=uep(F)\}
$$

\noindent forms an infinite sequence of random variables  

\begin{equation}\label{H2}
(N_{k,2})_{k\geq 1}=\{N_{1,2}, N_{2,2}, \cdots \},
\end{equation}

\bigskip
\noindent  Indeed, by denoting the event $uep(F)$ will not hit by any $X_j$, $j\geq 1$ by $A_n$, we have for $n=1$

\begin{eqnarray*}
\mathbb{P}(A_1)= \mathbb{P}\biggr(\bigcap_{j\geq 1} (X_j\neq uep(F))\biggr)=(1-p_2)^{+\infty}=0,
\end{eqnarray*}

\noindent and for any $n\geq 1$,

\begin{eqnarray*}
\mathbb{P}(A_n)&=& \sum_{k\geq 1} \mathbb{P}(A_n / (N_{n-1,2}=2)) \mathbb{P}(N_{n-1,2}=k)\\ 
&=&\mathbb{P}\biggr(\bigcap_{j\geq 1} (X_j\neq uep(F))\biggr)  \mathbb{P}(N_{n-1,2}=k)\\
&=& =(1-p_2)^{+\infty} \ \mathbb{P}(N_{n-1,2}=k)=0. \ \square
\end{eqnarray*}

\bigskip \noindent As well, if 

\begin{equation}\label{H1b}
\mathbb{P}\left(X=lep\left(F\right)\right)=p_1 \in \left] 0, 1\right[,
\end{equation}

\bigskip \noindent The sequence of hitting times of $lep(F)$ is infinite and is denoted as  

$$
(N_{k,1})_{k\geq 1}=\{N_{1,1}, N_{2,1}, \cdots \},
$$

\bigskip \noindent Throughout the paper, we suppose that $0<p_1+p_2<1$, otherwise $X$ would be a Bernoulli random variable and  the study would be quite simple. So, if both \eqref{H1a} and \eqref{H2} holds, we may define

$$
(N_{k,3})_{k\geq 1}=\{N_{1,1}, N_{2,1}, \cdots \},
$$

\bigskip
\noindent  as the random times in which neither endpoint is hit, that is

$$
\{j\geq 1, \ X_j\neq lep(F) \ and \ X_j\neq uep(F)\}.
$$

\bigskip \noindent We wish to describe the asymptotic theory of such sequences when they are defined. That theory reduces to studying sequences of negative binomial random variables. Let us make that first recall.

\begin{theorem} \label{theoHE01} Let $X$, $X_{1}$, $X_{2}$, \ldots  be \textit{iid} random variables with common \textit{cdf} $F$.\\ 
	
\noindent (a) If $uep(F) \in \mathbb{R}$ and \eqref{H1a} holds, the sequence $\left(N_{k,2}\right)_{k \geq 1}$ is an infinite sequence and for all $k \geq 1$, $N_{k,2}$ follows a negative binomial law of parameters and $k$  and $\overline{p}_2$
	
$$
N_{k,2} \sim \beta N \left(k, p_2\right).
$$
	
\bigskip\noindent (b) If  $lep(F) \in \mathbb{R}$ and \eqref{H1b} holds. Then the sequence $\left(N_{k,1} \right)_{k \geq 1}$ is infinite  and and for all $k \geq 1$, $N_{k,1}$ follows a negative binomial law of parameters and $k$  and $\overline{p}_1$
	
$$
N_{k,1} \sim \beta N \left(k, p_{1}\right).
$$
	
\bigskip \noindent (c) \eqref{H1a} and \eqref{H2} hold and $0 < p_1+p_2< 1 $, the sequence of $\left(N_{k,3}  \right)_{k \geq 1}$ is infinite and 
	
$$
\forall k \geq 1, N_{k,3} \sim \mathcal{\beta N} \left(k, p_{3}\right).
$$
\end{theorem}

\bigskip \noindent \textbf{Proof}. It is enough to prove this for one case. For example, if \eqref{H1a} holds, $Z_{1,2}=N_{1,2}$ follows a geometric law of parameter $p_{2}$, that is $Z_{1,2} \sim \mathcal{G}\left(p_{2}\right)$ and for all $h \geq 1$,

$$
\mathbb{P}\left(Z_{1,2}=h \right)= \overline{p}_{2}^{h-1} p_{2},
$$

\bigskip \noindent with 

$$
\overline{p}_{i}=1-p_{i}, \ i\in \{1, 2, 3 \}.
$$

\bigskip \noindent Once $Z_{1,2}$ is observed, the next hitting time is achieved after $Z_{2,2}$ trials and 

$$
N_{1,2}=Z_{1,2} + Z_{2,2}
$$

\bigskip \noindent with $Z_{1,2}$ independent of $Z_{2,2}$. By induction, we get $Z_{1,2}, Z_{2,2}, \cdots, Z_{k,2}, \cdots$ independent such that each 
$Z_{k,2} \sim \mathcal{G}\left(p_{2}\right)$ so that for $k \geq 1$,

$$
N_{k,2}=Z_{1,2}+ \cdots + Z_{k,2}.
$$

\bigskip \noindent This proves that $N_{k,2} \sim \beta N \left(k, p_{2}\right)$. $\square$

\bigskip \noindent From the previous proof, we saw that we have, for $i \in \{1, 2, 3 \}$, that there exists a sequence of independent and geometric random variables $\left(Z_{k,i} \right)_{k \geq 1}$ such that for all $k \geq 1$.

$$
N_{k,i}=Z_{1,i}+ \cdots + Z_{k,i}, \  the \\ Z_{\circ,i}'s \ \ are \ \ independent, \ Z_{\circ,i} \sim \mathcal{G}\left( p_{i} \right).
$$
  
\bigskip \noindent From there, the classical limit theorems for  sums of  \textit{iid} random variables apply as follows.\\ 

\noindent We denote, for $i \in \{1, 2, 3 \}$

$$
\nu_{i}=\mathbb{E}(Z_{1,i})=1/p_i,
$$

$$
\sigma_{i}^2=\mathbb{V}ar(Z_{1,i})=\overline{p}_{i}/(p_i^2),
$$

$$
\gamma_{i}=\mathbb{E}(Z_{1,i}-\mathbb{E}Z_{1,i})^3=\frac{\overline{p}_i(2-p_i)}{p^3_i}.
$$

\bigskip \noindent We have

\begin{theorem}\label{theo02} Upon \eqref{H1a} holds for $i=1$, \eqref{H1b} holds for $i=2$ and both \eqref{H1a} and \eqref{H1b} hold for $i=3$ with $p^{(3)} < 1$, we have  the following laws for $i \in \{ 1, 2, 3 \}$.\\

\noindent (a) \textbf{Central limit theorem (CLT)}.

$$
N_{k}^{(i, \star)}=\frac{N_{k}^{(i)}-k \nu_{i}}{\sigma_{i} \sqrt{k}} \rightsquigarrow \mathcal{N}\left(0,1\right)
$$

\bigskip \noindent (b) \textbf{Iterated logarithm law (LIL)}.

$$
\lim_{k \to + \infty} \frac{N_{k}^{(i)} - k \nu_{i}}{\sqrt{2k\sigma_{i}\log \log k}} =+1
$$

\bigskip \noindent and 

$$
\limsup_{k \to + \infty} \frac{N_{k}^{(i)} - k \nu_{i}}{\sqrt{2k\sigma_{i}\log \log k}} =+1 
$$

\bigskip \noindent (c) \textbf{Berry-Essen approximation (BEA)}.

$$
\sup_{x \in \mathbb{R}} \left| \mathbb{P}\left(N_{k}^{(i, \star)} \leq x \right)-\frac{1}{\sqrt{2 \pi}} \int_{- \infty}^{x} e^{-t^{2}/2} dt   \right| \leq 
\frac{36 \gamma_i}{\sqrt{k}}.
$$

\end{theorem}

\bigskip
\noindent The limit theorems (CLT), (LIL) and (BEA) are to be found in graduate probability textbooks like \cite{loeve}, \cite{chung}, \cite{gutt}, etc., in  particular in \cite{ips-mfpt-ang} (Theorem 20 page 237, Theorem 22 page 272, Theorem 21 page 253).\\

\bigskip \noindent If both endpoints are finite and atoms for the probability law, the \textit{cdf} may be represented as follows 
\begin{figure}[htbp]
	\caption{The two finite endpoints are both atoms}
	\centering
		\includegraphics[scale=0.75]{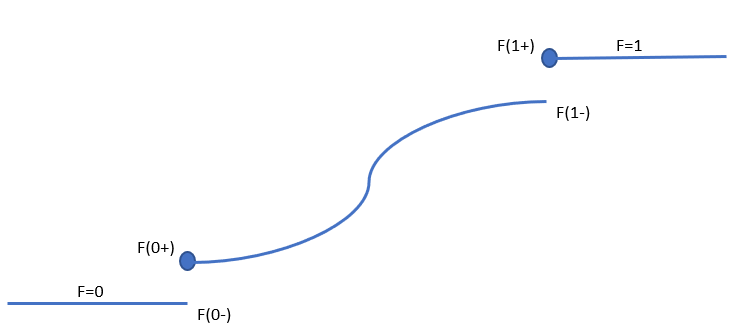}
	\label{fig1}
\end{figure}

\noindent In such cases, it would be interesting to compare the number of hitting times of one endpoint with the same number of the other endpoints or with number of hitting times of the non-endpoint zone.\\

\noindent Let us proceed by a multinomial approach. For each observation time $i \geq 1$, we simultaneously check which of the three events $A_{i}=\left(X_{i}=lep (F)\right)$, $B_{i}=\left(X_{i} \notin \{lep(F), uep(F) \}\right)$ or $C_{i}=\left(X_{i}=uep(F)  \right)$ occurs.

\noindent At time $n \geq 1$, $M_{n,i}$ is the number of occurrences of the $A_{i}$'s, of the $C_{i}$'s and finally of the $B_{i}$'s. So, we have that $M_{n,1}$ is number of hitting of $lep(F)$, $M_{n,2}$ the number of hitting of $uep(F)$ and $M_{n,3}$ the number of times neither $lep(F)$ nor $uep(F)$ are hit.\\

\noindent The vector $M_{n}= \left(N_{M,1}, M_{n,2}, M_{n,3} \right)$ follows a three dimensional multinomial law of parameters $\left(n,p\right)$, with $p=\left(p_{1}, p_{2}, p_{3}\right)$. We suppose that $p_{i}> 0$, $i \in \{ 1, 2, 3 \}$ and $p_{1}+p_{2}+p_{3}=1$.\\

\noindent We focus on compare $M_{n,1}$ and $M_{n,2}$ for large values of $n$, the results being extended to any pair of $\{(M_{n,1})_{n\geq 1}, (M_{n,2})_{n\geq 1}, (M_{n,3})_{n\geq 1}\}$.

\bigskip \noindent Let us use the classical result of the asymptotic law of $\left(M_{n} \right)_{n \geq 0}$.   
 
\begin{proposition} \label{rec-02} Given the context described above and the assumptions set above, we have 
	
$$
M_{n}^{\star}= \left(\frac{M_{n,1}-n p_{1}}{\sqrt{n p_{1}}}, \frac{M_{n,2}-n p_{2}}{\sqrt{n p_{2}}}, \frac{M_{n,3}-n p_{3}}{\sqrt{n p_{3}}} \right)^{t}
$$	  

\bigskip \noindent converges to a centered Gaussian vector $Z=\left(Z_{1}, Z_{2}, Z_{3}\right)^{t}$ with variance-covariance $\Sigma=\left( \sigma_{ij}  \right)_{1\leq i,j \leq 3}$:

$$
\left\{ 
\begin{array}{cc}
\sigma_{ii}=1-p_{i}&  \\
\sigma_{ij}=- \sqrt{p_{i}p_{j}},& i \neq j.
\end{array}
\right.
$$
\end{proposition}

\bigskip \noindent We infer the following laws. 

\begin{corollary} \label{rec-03} Given the notation above, we have 

\noindent (a)

$$
\sqrt{n} \left( \frac{M_{n,1}}{M_{n,2}}- \frac{p_{1}}{p_{2}}  \right) \sim \mathcal{N} \left(0, \gamma_{1,2}^{2}\right)
$$	

\bigskip \noindent with 

\begin{eqnarray*}
\gamma_{1,2}^{2}&=& (p_1/p_2)(p_2^2 \overline{p}_1+ p_1^2 \overline{p}_2 +2 (p_1p_2)^{3/2}).
\end{eqnarray*}

\bigskip \noindent (b)

$$
\sqrt{n} \left( \frac{M_{n,1}-M_{n,2}}{n}- \left( p_{1} -p_{2} \right)  \right) \sim \mathcal{N} \left(0, \delta_{1,2}^{2}\right)
$$ 

\bigskip \noindent with 

\begin{eqnarray*}
\delta_{1,2}^{2}=p_1 \overline{p}_1+ p_2 \overline{p}_2 +2 (p_1p_2).
\end{eqnarray*}

\end{corollary}

\bigskip \noindent Before we provide the proof of Corollary \ref{rec-03}, we give some of its important consequences. We are going to get the following. If the hitting probability of one of the endpoints is greater than the other by $\Delta p=|p_1-p_2|>0$, for any $\beta in ]0,\Delta p[$, its number of hitting times
is greater than the other counterpart of more that $n\beta$, that is (if $p_1>p_2$ for example)

$$
\liminf_{n \rightarrow +\infty} \mathbb{P}(M_{n,1}-M_{n,1}\geq n\beta)=100\%.
$$

\bigskip \noindent To see this, suppose that $\Delta=p_1-p_2>0$ and $0<\beta<\Delta$ and set  $I(n)=(M_{n,1}-M_{n,1}\geq n\beta)$, $n\geq 1$. By part (b) of Corollary \ref{rec-03}, and by $G$ the \textit{cdf} of absolute value the standard Gaussian random variable and by denoting

$$
C_n=\left| \sqrt{n} \left( \frac{M_{n,1}-M_{n,2}}{n}- \left( p_{1} -p_{2} \right)  \right)\right|.
$$ 

\noindent we have for any $t \in \mathbb{R}$, $\mathbb{P}(C_n\geq t)-(1-G(t))) \rightarrow 0$ and by a classical result in Weak convergence (when the \textit{limit cdf} is continuous, see for example \cite{ips-wcia-ang} [page 107, Chapter 4, Point (5)], we get

$$
\delta_n=\sup_{t \in \mathbb{R}} \left| \mathbb{P}(C_n\geq t)-(1-G(t)) \right| \rightarrow 0.
$$

\bigskip Let us apply this to 

\begin{eqnarray*}
\mathbb{P}(I(n))&=&\mathbb{P}\left(\frac{M_{n,1}-M_{n,2}}{n}\leq \Delta p\right)\\
&=&\mathbb{P}\left(\sqrt{n}\left(\frac{M_{n,1}-M_{n,2}}{n}-(p_1-p_2) \right)\leq \sqrt{n}(\beta -\Delta p)\right).
\end{eqnarray*}
 
\noindent We remark that for each $n\geq 1$, $-x_n=\sqrt{n}(\beta -\Delta p)$ is negative and $x_n \rightarrow +\infty$ as $n \rightarrow +\infty$. So we have

\begin{eqnarray*}
\mathbb{P}(I(n)^c)&\leq&\mathbb{P}\left(\left|\sqrt{n}\left(\frac{M_{n,1}-M_{n,2}}{n}-(p_1-p_2) \right)\right|\geq x_n\right)\\
&\leq& |\mathbb{P}(C_n\geq x_n)-(1-G(x_n))| + (1-G(x_n))\\
&\leq& \delta_n + (1-G(x_n)) \rightarrow 0,
\end{eqnarray*}

\noindent as $n \rightarrow +\infty$. We have proved that

$$
\liminf_{n \rightarrow +\infty} \mathbb{P}(I(n))=100\%. \ \square
$$

\bigskip \noindent \textbf{Proof of Corollary \ref{rec-03}}. We use the Skorohod-Wichura theorem (See \cite{skorohod} and \cite{wichura} and \cite{gslo-skorohodWichira-A} for a brief proof) and place ourselves on a probability space holding a sequence $\left(\overline{N}_{n}^{\star}\right)_{n \geq 1}$ and $\overline{Z}$ such that  

$$
\forall n \geq 1, \overline{N}_{n}^{\star}=_{d} M_{n}^{\star}, \ Z =_{d} \overline{Z}
$$

\bigskip \noindent and

$$
\overline{N}_{n}^{\star} \longrightarrow_{\mathbb{P}} \overline{Z}.
$$ 

\bigskip \noindent So, we may and do take $M_{n}^{\star}=\overline{N}_{n}^{\star}$ and $Z=\overline{Z}$. We have 

$$
\frac{M_{n,1}}{M_{n,2}}=\frac{\sqrt{np_{1}} \left( Z_{1}+o_{p}(1) \right)+ np_{1} }{\sqrt{np_{2}} \left( Z_{2}+o_{p}(1) \right)+ np_{2} }
$$

\bigskip \noindent From there, as proved in the Appendix,

\begin{eqnarray*}
\sqrt{n} \left( \frac{M_{n,1}}{M_{n,2}} - \frac{p_{1}}{p_{2}}\right)&=&\left(\frac{p_{1}}{p_{2}} \right)^{1/2} \left(p_{2}Z_{1} -p_{1}Z_{2}\right) \\
&\sim & \mathcal{N} \left(0, \gamma_{1,2}^{2} \right),
\end{eqnarray*}

\bigskip \noindent with 

\begin{eqnarray*}
\gamma_{1,2}^2&=&\mathbb{V}ar \left(\left(p_{1}p_{2}\right)^{1/2} \left(p_{2}Z_{1} - p_{1}Z_{2}\right)\right) \\
&=&\left(p_{1}/p_{2}  \right) \left[p_{2}^{2} \mathbb{E}Z_{1}^{2} +p_{1}^{2} \mathbb{E}Z_{2}^{2} - 2p_{1}p_{2} \mathbb{E} \left(Z_{1} Z_{2}  \right) \right].
\end{eqnarray*}

\bigskip \noindent Finally,

$$
\gamma_{1,2}^{2}=\left(p_{1} / p_{2}\right) \left[ p_{2}^{2}\left(1-p_{1}\right)+  p_{1}^{2}\left(1-p_{2}\right)+2\left(p_{1} p_{2}\right)^{3/2} \right].
$$
 
\bigskip \noindent We also have 
 
\begin{eqnarray*}
M_{n,1} - M_{n,2}&=&\left( \sqrt{n p_{1}}\left( Z_{1} + o_{p}(1) \right) +n p_{1} \right)- \left( \sqrt{n p_{2}}\left( Z_{2} + o_{p}(1) \right) +n p_{2} \right)\\
&=&n\left(p_{1}-p_{2}\right)+ \sqrt{n}\left( \sqrt{p_{1}}\left( Z_{1} + o_{p}(1) \right)  - \sqrt{p_{2}}\left( Z_{2} + o_{p}(1) \right)  \right).
\end{eqnarray*}
 
\bigskip \noindent So, 
 
$$
\sqrt{n}\left( \frac{M_{n,1}-M_{n,2}}{n} - \left(p_{1} - p_{2}\right) \right)= \sqrt{p_{1}}Z_{1} - \sqrt{p_{2}}Z_{2}+ o_{p}(1).
$$
 
\bigskip \noindent Since

\begin{eqnarray*}
\delta_{1,2}^{2}&=&\mathbb{V}ar \left( \sqrt{p_{1}}Z_{1} - \sqrt{p_{2}}Z_{2}  \right)\\
&=&p_{1}\left(1-p_{1}\right)+p_{2}\left(1-p_{2}\right)+2p_{1}p_{2}.
\end{eqnarray*}
  
\bigskip \noindent So,
  
$$
\sqrt{n}\left( \frac{M_{n,1}-M_{n,2}}{n} - \left(p_{1} - p_{2}\right) \right) \sim \mathcal{N}\left(0, \delta_{1,2}^{2}\right). \ \square
$$

\bigskip \noindent \textbf{Appendix}\\

\begin{eqnarray*}
\frac{M_{n,1}}{M_{n,2}} - \frac{p_{1}}{p_{2}}&=&\frac{\sqrt{n p_{1}}\left(M_{n}^{\star}(1) + np_{1} \right) }{ \sqrt{n p_{2}}\left(M_{n}^{\star}(2) + np_{2} \right)} - \frac{p_{1}}{p_{2}}\\
&=&\frac{p_{2} \sqrt{p_{1}p_{2}} \sqrt{n} M_{n}^{\star}(1) + np_{1}p_{2} - \left(\sqrt{n}p_{1}\sqrt{p_{1}}p_{2} M_{n}^{\star}(2)+np_{1}p_{2} \right)}{np_{2}\left(M_{n}^{\star}(2) / \left(\sqrt{np_{2}} \right)+1 \right) }\\
&=&\frac{\sqrt{p_{1}p_{2}} \left( p_{2} M_{n}^{\star} (1) - p_{1} M_{n}^{\star}(2)  \right)}{p_{2}\sqrt{n} \left(1 + M_{n}^{\star}(1)/\sqrt{n p_{2}}  \right) }\\
&=&\frac{\sqrt{p_{1}p_{2}} \left(p_{2}\left(Z_{1}+o_{P}(1) \right) - p_{1}\left(Z_{2}+o_{P}(1) \right) \right)}{p_{2}\sqrt{n} \left( 1 + \left(Z_{1}+o_{P}(1) \right) /\sqrt{n p_{2}}  \right)}.
\end{eqnarray*}

\bigskip \noindent Hence 
 
$$
\sqrt{n}\left( \frac{M_{n,1}}{M_{n,2}} - \frac{p_{1}}{p_{2}}\right)=\left(\frac{p_{1}}{p_{2}}\right)^{1/2} \left( p_{2}Z_{1} - p_{1}Z_{2}  \right) + o_{P}(1).
$$
 
\bigskip \noindent So,
 
$$
\sqrt{n}\left(\frac{M_{n,1}}{M_{n,2}} - \frac{p_{1}}{p_{2}}\right) \sim  \mathcal{N}\left(0, \gamma_{1,2}^{2}\right)
$$
 
\bigskip \noindent with

\begin{eqnarray*}
\gamma_{1,2}^{2}&=&\mathbb{V}ar\left( \left(p_{1} / p_{2} \right)^{1/2} \left(p_{2}Z_{1} - p_{1}Z_{2} \right) \right)\\
&=&\left(p_{1} / p_{2} \right)\left[p_{2}^{2}\mathbb{V}arZ_{1} + p_{1}^{2}\mathbb{V}arZ_{2} - 2p_{1}p_{2}\mathbb{C}ov\left( Z_{1}Z_{2} \right) \right]\\
&=&\left(p_{1} / p_{2} \right)\left[p_{2}^{2} \left(1 - p_{1}  \right) + p_{1}^{2}\left(1 - p_{2} \right) + 2\left( p_{1}p_{2} \right)^{3/2} \right].
\end{eqnarray*}

\section{A simple application} \label{sec3}

\noindent Let $X$ follows a binomial $\beta(r,\alpha)$, $r\geq 1$ and $\alpha \in ]0,1[$. The endpoints are $lep(F)=0$ and $uep(F)=r$ which are atoms of probability $(1-\alpha)^r$ and $\alpha^r$.\\

\noindent If $1-\alpha<\alpha$, i.e., $\alpha>1/2$, we have, by Theorem , for any $0<\delta<\alpha^r-(1-\alpha)^r$.

\begin{equation}
\liminf \mathbb{P}(M_{2,n}-M_{1,n}> n\delta)=100\%. \label{bcUepDominant}
\end{equation}

\bigskip
\noindent If $1-\alpha>\alpha$, i.e., $\alpha<1/2$, we have, by Part (b) of Corollary \ref{rec-03}, for any $0<\delta<(1-\alpha)^r-\alpha^r$

\begin{equation}
\liminf \mathbb{P}(M_{1,n}-M_{2,n}> n\delta)=100\%. \label{bcLepDominant}
\end{equation}

\bigskip \noindent If $\alpha=1/2$, $p_1=p_2=\alpha^3=(1-\alpha)^r=2^{-r}$. By applying Part (a) of Corollary \ref{rec-03}, we have

$$
\gamma_{1,2}=2 \gamma^2=2 2^{-2r}=2^{-2r+1}.
$$
 
\bigskip
\noindent and

$$
2^{r+1/2}\left(\left(M_{1,n}/M_{2,n}\right)-1\right) \rightsquigarrow \mathcal{N}(0,1).
$$

\bigskip \noindent By defining $z_{u}$ the critical points of a standard normal \textit{cdf}, that is $\mathbb{P}(\mathcal{N}(0,1)\leq z_u)=u$ for $u \in ]0,1[$, we get
confidence interval :

\begin{equation}
\forall u \in ]0,1[, \ \mathbb{P}\left(1-\frac{{2^{r}z_{1-u/2}}}{\sqrt{2n}} \leq \frac{M_{1,n}}{M_{2,n}} \leq 1+\frac{2^{r}z_{1-u/2}}{\sqrt{2n}}\right) \approx 1-u. 
\label{icequalityProbability}
\end{equation}

\end{document}